\documentclass{amsart}

\usepackage{tikz}
\usepackage{xcolor}
\usepackage{caption}

\newtheorem{theorem}{Theorem}[section]
\newtheorem{lemma}[theorem]{Lemma}

\begin{document}

\title{On the Morley's triangle}

\author{V.E.S. Szabó}
\address{Department of Analysis, Institute of
Mathematics, Budapest University of Technology and Economics, 
Műegyetem rkp. 3., Budapest, Hungary, H-1111}
\email{bmesszabo@gmail.com}


\date{August 21, 2022.}


\begin{abstract}
We give a complete investigation of Morley's trisector theorem. 
If the intersections of the half lines starting from the 
adjacent vertices of a triangle form an equilateral 
triangle for an arbitrary triangle, 
then the half lines are the angle trisectors.
To derive the result we use elementary trigonometry, Taylor series 
expansions, and solve systems of polynomial equations step by step.
As a byproduct we get a probably new trigonometric 
identity.
\end{abstract}

\maketitle

\section{Introduction}
\label{sec:Intro} 
The Morley's trisector theorem states 
that in any triangle, the three points of 
intersection of the adjacent angle trisectors 
form an equilateral triangle. 
More information can be found in \cite{COXETER1967Geometry}. 
We show that only the trisector configuration gives equilateral triangle.
To prove this we use trigonometry 
to derive a system of equations. Then applying 
the multivariable Taylor series expansion we 
obtain systems of polynomial equations. To find the 
solutions we perform direct calculations and use the 
resultant 
\cite{COX2015Ideals} 
as a tool.  
The remainder of paper is organized as follows. 
In Section \ref{sec:Trig} we derive a trigonometrical equation 
which plays key role in our calculations.  
In Section \ref{sec:Solution} applying 
multivariable Taylor series expansion we transform 
the problem into systems of polynomial equations. 
We solve the systems step by step.

\section{Trigonometrical equations}
\label{sec:Trig}
\noindent In the next lemma we note a consequence of 
the law of sines. If $A$ and $B$ are points, then denote  
$AB$ the length of side $AB$. 

\begin{lemma} \label{lemma:sinelaw}
Let $ABC_{\Delta}$ be an arbitrary triangle. Denote 
$ABC\angle:=\alpha$, $ACB\angle:=\beta$, and let $BC=1$. 

Then we have
    \begin{equation} \label{eq:sinelaw}
        AB=\frac{\sin(\beta)}{\sin(\alpha+\beta)},
        \quad AC=\frac{\sin(\alpha)}{\sin(\alpha+\beta)}.
    \end{equation}
\end{lemma}
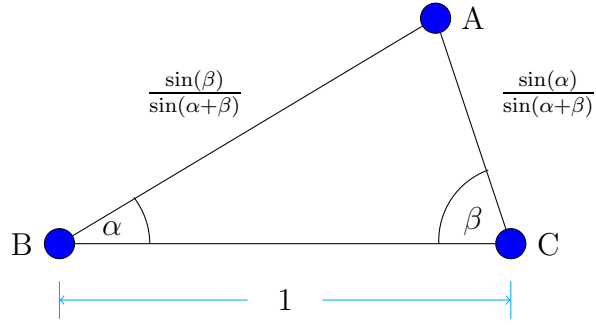
\begin{figure}[ht] 
\centering
\begin{tikzpicture}
  \draw (-4,0) -- (2,0);
  \draw (-4,0) -- (1,3);
  \draw (2,0) -- (1,3);
  \filldraw[fill=blue] (-4,0) circle [radius=2mm];
  \filldraw[fill=blue] (2,0) circle [radius=2mm];
  \filldraw[fill=blue] (1,3) circle [radius=2mm];
  
  \draw (-2.8,0) arc [start angle=0, end angle=38, radius=10mm];
  \draw (1.7,0.98) arc [start angle=110, end angle=182, radius=10mm];
  
  \draw (-4.5, 0) node {{\Large B}};
  \draw (2.5, 0) node {{\Large C}};
  \draw (1.5, 3) node {{\Large A}};
  \draw (-3.3, 0.2) node {{\Large $\alpha$}};
  \draw (1.5, 0.3) node {{\Large $\beta$}};
  
  \draw[<-, cyan] (-4, -0.75) -- (-1.5, -0.75); 
  \draw[cyan] (-4, -1) -- (-4, -0.5);
  \draw (-1, -0.75) node {{\Large $1$}};
  \draw[->, cyan] (-0.5, -0.75) -- (2, -0.75);
  \draw[cyan] (2, -1) -- (2, -0.5);
  
  \draw (-2.2, 2) node {{\Large $\frac{\sin(\beta)}{\sin(\alpha+\beta)}$}};
  \draw (2.5, 2) node {{\Large $\frac{\sin(\alpha)}{\sin(\alpha+\beta)}$}};
\end{tikzpicture}
\caption{Length of sides}
\end{figure}

\newpage
Now we derive the key equation.
\begin{lemma}\label{eq:keyeq}
    Let $EFH_{\Delta}$ be an arbitrary triangle. Denote $HEF\angle:=\alpha$,
$EFH\angle:=\beta$, $FHE\angle:=\gamma$. Let $G,I,J$ be inner points
of $EFH_{\Delta}$. Denote $GEF\angle:=t_{1}\alpha$, $GFE\angle:=t_{2}\beta$,
$IFH\angle:=t_{3}\beta$, $IHF\angle:=t_{4}\gamma$, $JHE\angle:=t_{5}\gamma$,
$JEH\angle:=t_{6}\alpha$, where $t_{i}>0$ $(i=1,\ldots,6)$, $t_{2}+t_{3}<1$,
$t_{4}+t_{5}<1$, $t_{6}+t_{1}<1$ and $t_{i}$ $(i=1,\ldots,6)$
are independent of $\alpha,\beta,\gamma$ and of the length of the
sides of $EFH_{\Delta}$.

Here $GI=IJ$ for every $0<\alpha,\beta$, $\alpha+\beta<\pi$ 
if and only if
\begin{equation}
A(\alpha,\beta)=0
\end{equation}
for every $0<\alpha,\beta$, $\alpha+\beta<\pi$, 
where
\begin{align*}
&A(\alpha,\beta):  =\\  
&\sin^{2}(t_{1}\alpha)\sin^{2}(\alpha+\beta)\sin^{2}(-t_{3}\beta-t_{4}\pi+t_{4}\alpha+t_{4}\beta) 
\sin^{2}(t_{5}(\alpha+\beta-\pi)-t_{6}\alpha)\\
 &   +\sin^{2}(\alpha)\sin^{2}(t_{4}(\pi-\alpha-\beta))\sin^{2}(t_{1}\alpha+t_{2}\beta) 
\sin^{2}(-t_{5}\pi+t_{5}\alpha+t_{5}\beta-t_{6}\alpha)\\
 &   +2\sin(t_{1}\alpha)\sin(\alpha)\sin(t_{4}(\pi-\alpha-\beta))\cos((-1+t_{2}+t_{3})\beta)
 \sin(t_{1}\alpha+t_{2}\beta) \\
 &\cdot    \sin(\alpha+\beta)\sin(-t_{3}\beta-t_{4}\pi+t_{4}\alpha+t_{4}\beta)
   \sin^{2}(-t_{5}\pi+t_{5}\alpha+t_{5}\beta-t_{6}\alpha)\\
 &  -\sin^{2}(\alpha)\sin^{2}(t_{3}\beta)\sin^{2}(t_{1}\alpha+t_{2}\beta)\sin^{2}(-t_{5}\pi+t_{5}\alpha+t_{5}\beta-t_{6}\alpha)\\
 &   -\sin^{2}(\beta)\sin^{2}(t_{6}\alpha)\sin^{2}(t_{1}\alpha+t_{2}\beta)\sin^{2}(-t_{3}\beta-t_{4}\pi+t_{4}\alpha+t_{4}\beta)\\
 &   +2\sin(\alpha)\sin(t_{3}\beta)\sin(\beta)\sin(t_{6}\alpha)\cos((-1+t_{4}+t_{5})(\pi-\alpha-\beta))\\
&\cdot \sin^{2}(t_{1}\alpha+t_{2}\beta)
 \sin(-t_{3}\beta-t_{4}\pi+t_{4}\alpha+t_{4}\beta)
 \sin(-t_{5}\pi+t_{5}\alpha+t_{5}\beta-t_{6}\alpha).
\end{align*}
\end{lemma}

\begin{figure}[ht]
\centering
\begin{tikzpicture}
  \draw (-8,0) -- (3,0);
  \draw (-8,0) -- (0,7);
  \draw (0,7) -- (3,0);
  
  \draw[magenta] (-8,0) -- (-1.5,2); 
  \draw[green] (-1.5,2) -- (3,0);  
  
  \draw[magenta] (-8,0) -- (-2,3); 
  \draw[orange] (-2,3) -- (0,7);
  
  \draw[green] (3,0) -- (-1,2.5); 
  \draw[orange] (-1,2.5) -- (0,7);
  
   \draw[magenta] (-2,3) -- (-1.5,2);
  \draw[green] (-1.5,2) -- (-1,2.5);
  \draw[orange] (-1,2.5) -- (-2,3);
  
  \filldraw[fill=blue] (-8,0) circle [radius=2mm]; 
  \filldraw[fill=blue] (3,0) circle [radius=2mm]; 
  \filldraw[fill=blue] (0,7) circle [radius=2mm]; 
  \filldraw[fill=blue] (-1.5,2) circle [radius=2mm]; 
  \filldraw[fill=blue] (-2,3) circle [radius=2mm]; 
  \filldraw[fill=blue] (-1,2.5) circle [radius=2mm]; 

  \draw (-8, -0.5) node {{\Large E}};
  \draw (3, -0.5) node {{\Large F}};
  \draw (0.5, 7) node {{\Large H}};
  \draw (-1.5, 1.5) node {{\Large G}}; 
  \draw (-2.5, 3) node {{\Large J}}; 
  \draw (-0.5, 2.5) node {{\Large I}}; 
  
  \draw (-5.5, 0.3) node {{\Large $t_1\alpha$}};
  \draw (-5.75, 1.5)  node {{\Large $t_6\alpha$}};
  \draw (1.1, 0.3) node {{\Large $t_2\beta$}};
  \draw (2, 1.1) node {{\Large $t_3\beta$}};
  \draw (0.1, 5.7) node {{\Large $t_4\gamma$}};
  \draw (-1.1, 5.5) node {{\Large $t_5\gamma$}};
  
  \draw (-1, -0.5) node {{\Large $1$}};
  \draw[magenta] (-3.2, 0.8) node {{\Large 
    $\frac{\sin(t_2\beta)}{\sin(t_1\alpha+t_2\beta)}$}};
  \draw[green] (-0.3, 0.8) node {{\Large 
    $\frac{\sin(t_1\alpha)}{\sin(t_1\alpha+t_2\beta)}$}};
  \draw[green] (1.9, 2.3) node {{\Large 
    $\frac{\sin(\alpha)}{\sin(\alpha+\beta)}\frac{\sin(t_4\gamma)}{\sin(t_3\beta+t_4\gamma)}$}};
  \draw[orange] (1.6, 4.8) node {{\Large 
    $\frac{\sin(\alpha)}{\sin(\alpha+\beta)}\frac{\sin(t_3\beta)}{\sin(t_3\beta+t_4\gamma)}$}};
  \draw[orange] (-3.1, 4.8) node {{\Large 
    $\frac{\sin(\beta)}{\sin(\alpha+\beta)}\frac{\sin(t_6\alpha)}{\sin(t_5\gamma+t_6\alpha)}$}};
  \draw[magenta] (-5.6, 2.5) node {{\Large 
    $\frac{\sin(\beta)}{\sin(\alpha+\beta)}\frac{\sin(t_5\gamma)}{\sin(t_5\gamma+t_6\alpha)}$}};
  
  \draw (-4.8, 3.8) node {{\Large $\frac{\sin(\beta)}{\sin(\alpha+\beta)}$}};
  \draw (2.5, 3.3) node {{\Large $\frac{\sin(\alpha)}{\sin(\alpha+\beta)}$}};
\end{tikzpicture}
\caption{Length of sides}
\end{figure}
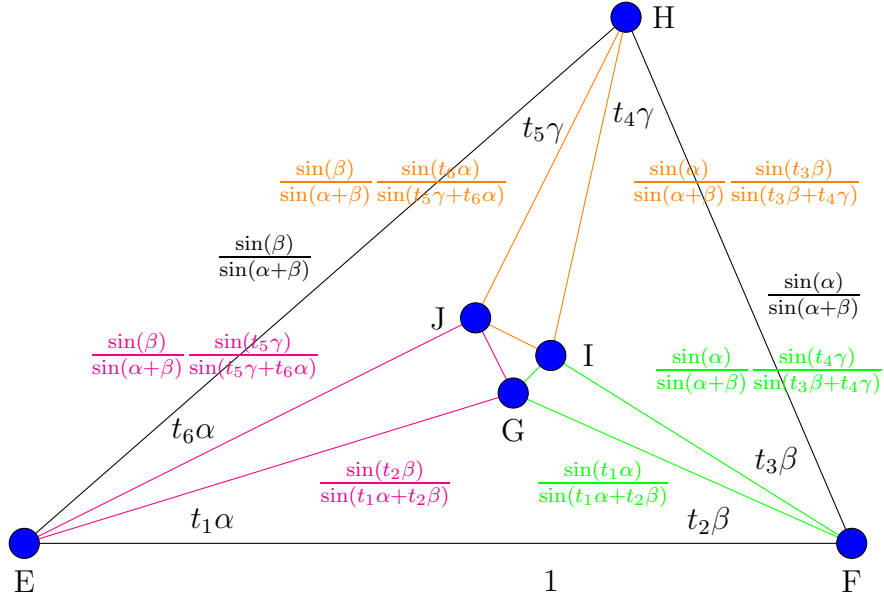

\newpage
\begin{proof}
Without the loss of generality we may assume that $EF:=1$. By
the law of sines we obtain
\begin{equation*}
FH=\frac{\sin(\alpha)}{\sin(\alpha+\beta)},
\end{equation*}
\begin{equation*}
HE=\frac{\sin(\beta)}{\sin(\alpha+\beta)}.
\end{equation*}
Applying the law of sines repeatedly it follows
\begin{equation*}
EG=\frac{\sin(t_{2}\beta)}{\sin(t_{1}\alpha+t_{2}\beta)},
\end{equation*}
\begin{equation*}
GF=\frac{\sin(t_{1}\alpha)}{\sin(t_{1}\alpha+t_{2}\beta)},
\end{equation*}
\begin{equation*}
FI=\frac{\sin(\alpha)}{\sin(\alpha+\beta)}\cdot\frac{\sin(t_{4}\gamma)}{\sin(t_{3}\beta+t_{4}\gamma)},
\end{equation*}
\begin{equation*}
IH=\frac{\sin(\alpha)}{\sin(\alpha+\beta)}\cdot\frac{\sin(t_{3}\beta)}{\sin(t_{3}\beta+t_{4}\gamma)},
\end{equation*}
\begin{equation*}
HJ=\frac{\sin(\beta)}{\sin(\alpha+\beta)}\cdot\frac{\sin(t_{6}\alpha)}{\sin(t_{5}\gamma+t_{6}\alpha)},
\end{equation*}
\begin{equation*}
JE=\frac{\sin(\beta)}{\sin(\alpha+\beta)}\cdot\frac{\sin(t_{5}\gamma)}{\sin(t_{5}\gamma+t_{6}\alpha)}.
\end{equation*}
Applying the law of cosine we get
\begin{align*}
GI^{2}  = & \left(\frac{\sin(t_{1}\alpha)}{\sin(t_{1}\alpha+t_{2}\beta)}\right)^{2}+\left(\frac{\sin(\alpha)}{\sin(\alpha+\beta)}\cdot\frac{\sin(t_{4}(\pi-\alpha-\beta))}{\sin(t_{3}\beta+t_{4}(\pi-\alpha-\beta))}\right)^{2}\\
   & -2\frac{\sin(t_{1}\alpha)}{\sin(t_{1}\alpha+t_{2}\beta)}\cdot\frac{\sin(\alpha)}{\sin(\alpha+\beta)}\cdot\frac{\sin(t_{4}(\pi-\alpha-\beta))}{\sin(t_{3}\beta+t_{4}(\pi-\alpha-\beta))} \\
   &\cdot\cos((1-t_{2}-t_{3})\beta)
\end{align*}
and
\begin{align*}
IJ^{2}  = & \left(\frac{\sin(\alpha)}{\sin(\alpha+\beta)}\cdot\frac{\sin(t_{3}\beta)}{\sin(t_{3}\beta+t_{4}(\pi-\alpha-\beta))}\right)^{2}\\
   & +\left(\frac{\sin(\beta)}{\sin(\alpha+\beta)}\cdot\frac{\sin(t_{6}\alpha)}{\sin(t_{5}(\pi-\alpha-\beta)+t_{6}\alpha)}\right)^{2}\\
   & -2\frac{\sin(\alpha)}{\sin(\alpha+\beta)}\cdot\frac{\sin(t_{3}\beta)}{\sin(t_{3}\beta+t_{4}(\pi-\alpha-\beta))}\cdot\frac{\sin(\beta)}{\sin(\alpha+\beta)}\\
   &\cdot\frac{\sin(t_{6}\alpha)}{\sin(t_{5}(\pi-\alpha-\beta)+t_{6}\alpha)}\cdot\cos((1-t_{4}-t_{5})(\pi-\alpha-\beta)).
\end{align*}
If $GI=IJ$ then $GI^2=IJ^2$, which yields
\begin{equation*}
A(\alpha,\beta)=0
\end{equation*}
for every $0<\alpha,\beta$, $\alpha+\beta<\pi$.
\end{proof}
We note that since $A(\alpha,\beta)$ is an analytical function we 
can drop the restrictions for $\alpha,\beta$ in the following. 

\section{Solution of the system of equations}
\label{sec:Solution}
\begin{theorem} 
If the intersections of the half lines starting from the 
adjacent vertices of a triangle form an equilateral 
triangle for an arbitrary triangle, 
then the half lines are the angle trisectors.    
\end{theorem}
\begin{proof}
Calculating the coefficient of $\alpha^{2}\beta^{2}$ in the 
Taylor series expansion of $A(\alpha,\beta)$ we arrive at the 
equation
\begin{equation*}
(-1+\cos^{2}(t_{4}\pi))(-1+\cos^{2}(t_{5}\pi))(t_{1}-t_{2})^{2}=0.
\end{equation*}
 From this we obtain
\begin{equation}
t_{1}=t_{2}.\label{eq:t12}
\end{equation}
In the same way, $IJ=JG$ implies
\begin{equation}
t_{3}=t_{4},\label{eq:t34}
\end{equation}
 and $JG=GI$ implies
\begin{equation}
t_{5}=t_{6}.\label{eq:t56}
\end{equation}
 Substituting (\ref{eq:t12}), (\ref{eq:t34}) and (\ref{eq:t56})
into $A(\alpha,\beta)$ we have
\begin{eqnarray*}
A(\alpha,\beta) & = & \sin^{2}(t_{1}\alpha)\sin^{2}(\alpha+\beta)\sin^{2}(-t_{3}\pi+t_{3}\alpha)\sin^{2}(-t_{5}\pi+t_{5}\beta)\\
 &  & +\sin^{2}(\alpha)\sin^{2}(t_{3}(\pi-\alpha-\beta))\sin^{2}(t_{1}\alpha+t_{1}\beta)\sin^{2}(-t_{5}\pi+t_{5}\beta)\\
 &  & +2\sin(t_{1}\alpha)\sin(\alpha)\sin(t_{3}(\pi-\alpha-\beta))\cos((-1+t_{1}+t_{3})\beta)\\
 &  & \cdot\sin(t_{1}\alpha+t_{1}\beta)\sin(\alpha+\beta)\sin(-t_{3}\pi+t_{3}\alpha)\sin^{2}(-t_{5}\pi+t_{5}\beta)\\
 &  & -\sin^{2}(\alpha)\sin^{2}(t_{3}\beta)\sin^{2}(t_{1}\alpha+t_{1}\beta)\sin^{2}(-t_{5}\pi+t_{5}\beta)\\
 &  & -\sin^{2}(\beta)\sin^{2}(t_{5}\alpha)\sin^{2}(t_{1}\alpha+t_{1}\beta)\sin^{2}(-t_{3}\pi+t_{3}\alpha)\\
 &  & +2\sin(\alpha)\sin(t_{3}\beta)\sin(\beta)\sin(t_{5}\alpha) 
 \cos((-1+t_{3}+t_{5})(\pi-\alpha-\beta))\\
 &  & \cdot\sin^{2}(t_{1}\alpha+t_{1}\beta)\sin(-t_{3}\pi+t_{3}\alpha)\sin(-t_{5}\pi+t_{5}\beta).
\end{eqnarray*}
 The coefficient of $\alpha^{2}\beta^{2}(\alpha+\beta)^{2}$ in the
Taylor expansion of $A(\alpha,\beta)$ about $(0,0)$ must be $0$,
hence   
\begin{align*}
& 2t_{3}t_{5}\cos((-1+t_{3}+t_{5})\pi)\sin(t_{3}\pi)\sin(t_{5}\pi)\\
&+\left((-1+t_{1}+t_{3})^{2}-t_{3}^{2}\right)\sin^{2}(t_{3}\pi)\sin^{2}(t_{5}\pi)-t_{5}^{2}\sin^{2}(t_{3}\pi)=0.
\end{align*}
If $(t_{1}^{*},t_{3}^{*},t_{5}^{*})$ is a solution of this equation,
then any permutation of it again a solution, because $t_{1},t_{3},t_{5}$
are independent of $\alpha,\beta,\gamma$. Denote $c_{i}:=\cos(t_{i}\pi)$,
$s_{i}:=\sin(t_{i}\pi)$, $(i=1,3,5)$. Thus we can write the following
system of equations
\begin{eqnarray}
-2t_{3}t_{5}c_{3}c_{5}s_{5}+((1-2t_{1}-2t_{3}+t_{1}^{2}+2t_{1}t_{3}+2t_{3}t_{5})s_{5}^{2}-t_{5}^{2})s_{3} & = & 0,\label{eq:Eeq1}\\
-2t_{5}t_{3}c_{5}c_{3}s_{3}+((1-2t_{1}-2t_{5}+t_{1}^{2}+2t_{1}t_{5}+2t_{5}t_{3})s_{3}^{2}-t_{3}^{2})s_{5} & = & 0,\label{eq:Eeq2}\\
-2t_{1}t_{5}c_{1}c_{5}s_{5}+((1-2t_{3}-2t_{1}+t_{3}^{2}+2t_{3}t_{1}+2t_{1}t_{5})s_{5}^{2}-t_{5}^{2})s_{1} & = & 0,\label{eq:Eeq3}\\
-2t_{5}t_{1}c_{5}c_{1}s_{1}+((1-2t_{3}-2t_{5}+t_{3}^{2}+2t_{3}t_{5}+2t_{5}t_{1})s_{1}^{2}-t_{1}^{2})s_{5} & = & 0,\label{eq:Eeq4}\\
-2t_{1}t_{3}c_{1}c_{3}s_{3}+((1-2t_{5}-2t_{1}+t_{5}^{2}+2t_{5}t_{1}+2t_{1}t_{3})s_{3}^{2}-t_{3}^{2})s_{1} & = & 0,\label{eq:Eeq5}\\
-2t_{3}t_{1}c_{3}c_{1}s_{1}+((1-2t_{5}-2t_{3}+t_{5}^{2}+2t_{5}t_{3}+2t_{3}t_{1})s_{1}^{2}-t_{1}^{2})s_{3} & = & 0.\label{eq:Eeq6}
\end{eqnarray}
Denote $S_{1}:=s_{1}^{2}$, $S_{3}:=s_{3}^{2}$, $S_{5}:=s_{5}^{2}$.
Multiplying (\ref{eq:Eeq1}) by $s_{3}$, (\ref{eq:Eeq2}) by $s_{5}$,
and equating the left hand sides, after some simplifications we obtain
\begin{equation*}
(-2t_{3}+2t_{1}t_{3})S_{5}S_{3}-t_{5}^{2}S_{3}=(-2t_{5}+2t_{1}t_{5})S_{3}S_{5}-t_{3}^{2}S_{5}.
\end{equation*}
 From this we get
\begin{equation*}
2(t_{5}-t_{3})(1-t_{1})S_{3}S_{5}=t_{5}^{2}S_{3}-t_{3}^{2}S_{5},
\end{equation*}
 which gives
\begin{equation}
[2(t_{5}-t_{3})(1-t_{1})S_{5}-t_{5}^{2}]S_{3}=-t_{3}^{2}S_{5},
\end{equation}
 and
\begin{equation}
S_{3}=-\frac{t_{3}^{2}S_{5}}{2(t_{5}-t_{3})(1-t_{1})S_{5}-t_{5}^{2}}.\label{eq:S3}
\end{equation}
 Similarly, we obtain from (\ref{eq:Eeq3}), (\ref{eq:Eeq4})
\begin{equation}
S_{1}=-\frac{t_{1}^{2}S_{5}}{2(t_{5}-t_{1})(1-t_{3})S_{5}-t_{5}^{2}}.\label{eq:S1}
\end{equation} 

Now consider the subsystem (\ref{eq:Eeq1}), (\ref{eq:Eeq4}), (\ref{eq:Eeq5}).
The assumptions imply that the coefficients of $s_{3},s_{5},s_{1}$
can not be zero. Indeed, assume that one of them is zero, say
\begin{equation}
(1-2t_{1}-2t_{3}+t_{1}^{2}+2t_{1}t_{3}+2t_{3}t_{5})s_{5}^{2}-t_{5}^{2}=0.\label{eq:Dis1}
\end{equation}
 Then from (\ref{eq:Eeq1}) we obtain $c_{3}c_{5}=0$. Assume that
$c_{3}=0$. Then we get
\begin{equation}
t_{3}=\frac{1}{2}.\label{eq:Dis2}
\end{equation}
Calculating the coefficient of $\alpha^{5}\beta^{2}$ in the Taylor
series expansion of $A(\alpha,\beta)$ we obtain
\begin{align}
&\frac{3}{2}t_{5}s_{3}s_{5}(t_{3}+t_{5}-1)(s_{3}c_{5}+c_{3}s_{5})\nonumber\\
&+\left(\frac{3}{2}t_{3}t_{5}(s_{3}s_{5}-c_{3}c_{5})s_{5}+s_{3}
((t_{1}-1)(t_{1}-2+3t_{3})s_{5}^{2}-\frac{3}{2}t_{5}^{2}
)\right)c_{3}=0.\label{eq:Eeq7}
\end{align}
 Substituting $t_{3}=1/2$ into (\ref{eq:Eeq7}) it follows
\begin{equation*}
\left(t_{5}-\frac{1}{2}\right)c_{5}=0.
\end{equation*}
 Hence, $t_{5}=1/2$. But this is a contradiction, because $1>t_{4}+t_{5}=t_{3}+t_{5}=1$
is impossible.

Assume that $c_{5}=0$, that is, $t_{5}=1/2$. As we noted earlier
if we have an equation for $t_{1},t_{3},t_{5}$ than any permutation
of variables again a valid equation. Swapping $3$ and $5$ in (\ref{eq:Eeq7})
we obtain
\begin{align}
&\frac{3}{2}t_{3}s_{5}s_{3}(t_{3}+t_{5}-1)(s_{3}c_{5}+c_{3}s_{5})\nonumber\\
&+\left(\frac{3}{2}t_{3}t_{5}(s_{3}s_{5}-c_{3}c_{5})s_{3}+s_{5}((t_{1}-1)(t_{1}-2+3t_{5})s_{3}^{2}-\frac{3}{2}t_{3}^{2})\right)c_{5}=0.\label{eq:Dis3}
\end{align}
Substituting $t_{5}=1/2$ into the above equation it implies
\begin{equation*}
\left(t_{3}-\frac{1}{2}\right)c_{3}=0,
\end{equation*}
 which is again a contradiction.

Returning to the system (\ref{eq:Eeq1}), (\ref{eq:Eeq4}), (\ref{eq:Eeq5})
we obtain
\begin{eqnarray}
c_{3}c_{5} & = & \frac{((1-2t_{1}-2t_{3}+t_{1}^{2}+2t_{1}t_{3}+2t_{3}t_{5})s_{5}^{2}-t_{5}^{2})s_{3}}{2t_{3}t_{5}s_{5}},\label{eq:E1}\\
c_{5}c_{1} & = & \frac{((1-2t_{3}-2t_{5}+t_{3}^{2}+2t_{3}t_{5}+2t_{5}t_{1})s_{1}^{2}-t_{1}^{2})s_{5}}{2t_{5}t_{1}s_{1}},\label{eq:E4}\\
c_{1}c_{3} & = & \frac{((1-2t_{5}-2t_{1}+t_{5}^{2}+2t_{5}t_{1}+2t_{1}t_{3})s_{3}^{2}-t_{3}^{2})s_{1}}{2t_{3}t_{1}s_{3}}.\label{eq:E5}
\end{eqnarray}
 Taking the product of equations we obtain
\begin{eqnarray}
c_{1}^{2}c_{3}^{2}c_{5}^{2} & = & \frac{1}{8t_{1}^{2}t_{3}^{2}t_{5}^{2}}\cdot((1-2t_{1}-2t_{3}+t_{1}^{2}+2t_{1}t_{3}+2t_{3}t_{5})s_{5}^{2}-t_{5}^{2})\nonumber \\
 &  & \times((1-2t_{3}-2t_{5}+t_{3}^{2}+2t_{3}t_{5}+2t_{5}t_{1})s_{1}^{2}-t_{1}^{2})\label{eq:product145}\\
 &  & \times((1-2t_{5}-2t_{1}+t_{5}^{2}+2t_{5}t_{1}+2t_{1}t_{3})s_{3}^{2}-t_{3}^{2}).\nonumber 
\end{eqnarray}
Taking the squares of equations (\ref{eq:E1}), (\ref{eq:E4}), (\ref{eq:E5})
we get
\begin{eqnarray}
c_{3}^{2}c_{5}^{2} & = & \frac{((1-2t_{1}-2t_{3}+t_{1}^{2}+2t_{1}t_{3}+2t_{3}t_{5})s_{5}^{2}-t_{5}^{2})^{2}s_{3}^{2}}{4t_{3}^{2}t_{5}^{2}s_{5}^{2}},\label{eq:E1sq}\\
c_{5}^{2}c_{1}^{2} & = & \frac{((1-2t_{3}-2t_{5}+t_{3}^{2}+2t_{3}t_{5}+2t_{5}t_{1})s_{1}^{2}-t_{1}^{2})^{2}s_{5}^{2}}{4t_{5}^{2}t_{1}^{2}s_{1}^{2}},\label{eq:E4sq}\\
c_{1}^{2}c_{3}^{2} & = & \frac{((1-2t_{5}-2t_{1}+t_{5}^{2}+2t_{5}t_{1}+2t_{1}t_{3})s_{3}^{2}-t_{3}^{2})^{2}s_{1}^{2}}{4t_{3}^{2}t_{1}^{2}s_{3}^{2}}.\label{eq:E5sq}
\end{eqnarray}
From (\ref{eq:product145}) and (\ref{eq:E1sq}), (\ref{eq:E4sq}),
(\ref{eq:E5sq}) we have
\begin{eqnarray}
c_{1}^{2}=& s_{5}^{2}((1-2t_{3}-2t_{5}+t_{3}^{2}+2t_{3}t_{5}+2t_{5}t_{1})s_{1}^{2}-t_{1}^{2})\nonumber \\ &\times\dfrac{((1-2t_{5}-2t_{1}+t_{5}^{2}+2t_{5}t_{1}+2t_{1}t_{3})s_{3}^{2}-t_{3}^{2})}{2t_{1}^{2}s_{3}^{2}((1-2t_{1}-2t_{3}+t_{1}^{2}+2t_{1}t_{3}+2t_{3}t_{5})s_{5}^{2}-t_{5}^{2})},\label{eq:c1sq1} 
\end{eqnarray}
\begin{eqnarray}
c_{3}^{2}=& s_{1}^{2}((1-2t_{1}-2t_{3}+t_{1}^{2}+2t_{1}t_{3}+2t_{3}t_{5})s_{5}^{2}-t_{5}^{2}) \nonumber \\   &\times\dfrac{((1-2t_{5}-2t_{1}+t_{5}^{2}+2t_{5}t_{1}+2t_{1}t_{3})s_{3}^{2}-t_{3}^{2})}{2t_{3}^{2}s_{5}^{2}((1-2t_{3}-2t_{5}+t_{3}^{2}+2t_{3}t_{5}+2t_{5}t_{1})s_{1}^{2}-t_{1}^{2})},\label{eq:c3sq1}
\end{eqnarray}
\begin{eqnarray}
c_{5}^{2}=& s_{3}^{2}((1-2t_{1}-2t_{3}+t_{1}^{2}+2t_{1}t_{3}+2t_{3}t_{5})s_{5}^{2}-t_{5}^{2})\nonumber \\  &\times\dfrac{((1-2t_{3}-2t_{5}+t_{3}^{2}+2t_{3}t_{5}+2t_{5}t_{1})s_{1}^{2}-t_{1}^{2})}{2t_{5}^{2}s_{1}^{2}((1-2t_{5}-2t_{1}+t_{5}^{2}+2t_{5}t_{1}+2t_{1}t_{3})s_{3}^{2}-t_{3}^{2})}.\label{eq:c5sq1}
\end{eqnarray}   
Denote $C_{1}:=c_{1}^{2}$, $C_{3}:=c_{3}^{2}$, $C_{5}:=c_{5}^{2}$.
Substituting (\ref{eq:S1}), (\ref{eq:S3}) into the above three expressions
we obtain
\begin{eqnarray}
C_{1}  = & -\dfrac{((1+t_{5}^{2}+4t_{1}t_{3}-2t_{3}-2t_{1})S_{5}-t_{5}^{2})}
{2((1-2t_{1}+t_{1}^{2}+2t_{1}t_{3}+2t_{3}t_{5}-2t_{3})S_{5}-t_{5}^{2})}\nonumber\\
&\times\dfrac{((1-2t_{3}+2t_{5}t_{1}-2t_{1}+2t_{1}t_{3}+t_{3}^{2})S_{5}-t_{5}^{2})}{((2t_{1}t_{3}+2t_{5}-2t_{3}t_{5}-2t_{1})S_{5}-t_{5}^{2})},\label{eq:C1}
\end{eqnarray}
\begin{eqnarray}
C_{3}  = & -\dfrac{((1+t_{5}^{2}+4t_{1}t_{3}-2t_{3}-2t_{1})S_{5}-t_{5}^{2})}
{2((1-2t_{3}+2t_{5}t_{1}-2t_{1}+2t_{1}t_{3}+t_{3}^{2})S_{5}-t_{5}^{2})}\nonumber\\
&\times\dfrac{((1-2t_{1}+t_{1}^{2}+2t_{1}t_{3}+2t_{3}t_{5}-2t_{3})S_{5}-t_{5}^{2})}{((2t_{1}t_{3}+2t_{5}-2t_{1}t_{5}-2t_{3})S_{5}-t_{5}^{2})},\label{eq:C3}
\end{eqnarray}
\begin{eqnarray}
C_{5}  = & \dfrac{((1-2t_{3}+2t_{5}t_{1}-2t_{1}+2t_{1}t_{3}+t_{3}^{2})S_{5}-t_{5}^{2})}
{2((1+t_{5}^{2}+4t_{1}t_{3}-2t_{3}-2t_{1})S_{5}-t_{5}^{2})t_{5}^{2}}\nonumber\\
&\times ((1-2t_{1}+t_{1}^{2}+2t_{1}t_{3}+2t_{3}t_{5}-2t_{3})S_{5}-t_{5}^{2}).\label{eq:C5}
\end{eqnarray}
Our next task is to determine $S_{5}$. From (\ref{eq:Eeq7}) it
follows
\begin{eqnarray}
&-\dfrac{3}{2}t_{3}s_{5}s_{3}(t_{3}+t_{5}-1)(s_{3}c_{5}+c_{3}s_{5})\nonumber\\
&=\left(\dfrac{3}{2}t_{3}t_{5}(s_{3}s_{5}-c_{3}c_{5})s_{3}+s_{5}\left((t_{1}-1)(t_{1}-2+3t_{5})s_{3}^{2}-\dfrac{3}{2}t_{3}^{2}\right)\right)c_{5}.\label{eq:S5eq1}
\end{eqnarray}
Using (\ref{eq:E1}) and (\ref{eq:S3}) we have
\begin{eqnarray}
(s_{3}s_{5}-c_{3}c_{5})s_{3} & = & S_{3}s_{5}-c_{3}c_{5}s_{3}\nonumber \\
 & = & S_{3}s_{5}-\frac{((1-2t_{1}-2t_{3}+t_{1}^{2}+2t_{1}t_{3}+2t_{3}t_{5})s_{5}^{2}-t_{5}^{2})S_{3}}{2t_{3}t_{5}s_{5}}\nonumber \\
 & = & S_{3}s_{5}+\frac{((1-2t_{1}-2t_{3}+t_{1}^{2}+2t_{1}t_{3}+2t_{3}t_{5})S_{5}-t_{5}^{2})t_{3}s_{5}}{2t_{5}\left(2(t_{5}-t_{3})(1-t_{1})S_{5}-t_{5}^{2}\right)}\nonumber \\
 & = & \frac{s_{5}t_{3}\left((1-t_{1})(1-t_{1}-2t_{3})S_{5}-t_{5}^{2}\right)}{2t_{5}\left(2(t_{5}-t_{3})(1-t_{1})S_{5}-t_{5}^{2}\right)}.\label{eq:S5eq1right1}
\end{eqnarray}
Using (\ref{eq:S3}) we have
\begin{equation}
(t_{1}-1)(t_{1}-2+3t_{5})s_{3}^{2}-\frac{3}{2}t_{3}^{2}=-\frac{t_{3}^{2}\left(2(1-t_{1})(2-t_{1}-3t_{3})S_{5}-3t_{5}^{2}\right)}{4(t_{5}-t_{3})(1-t_{1})S_{5}-2t_{5}^{2}}.\label{eq:S5eq1right2}
\end{equation}
From (\ref{eq:S5eq1right1}) and (\ref{eq:S5eq1right2}) we get
\begin{eqnarray}
&\left(\dfrac{3}{2}t_{3}t_{5}(s_{3}s_{5}-c_{3}c_{5})s_{3}+s_{5}\left((t_{1}-1)(t_{1}-2+3t_{5})s_{3}^{2}-\dfrac{3}{2}t_{3}^{2}\right)\right)c_{5}
\nonumber\\
&=\dfrac{c_{5}s_{5}t_{3}^{2}\left((1-t_{1})(t_{1}-5+6t_{3})S_{5}+3t_{5}^{2}\right)}{8(t_{5}-t_{3})(1-t_{1})S_{5}-4t_{5}^{2}}.\label{eq:S5eq1right}
\end{eqnarray} 
Substituting (\ref{eq:S5eq1right}) into (\ref{eq:S5eq1}) and dividing
by $c_{5}s_{5}t_{3}\neq0$ we obtain
\begin{equation}
-\frac{3}{2}(t_{3}+t_{5}-1)\left(S_{3}+s_{3}s_{5}\frac{c_{3}}{c_{5}}\right)=\frac{t_{3}\left((1-t_{1})(t_{1}-5+6t_{3})S_{5}+3t_{5}^{2}\right)}{8(t_{5}-t_{3})(1-t_{1})S_{5}-4t_{5}^{2}}.\label{eq:S5eq2}
\end{equation}
By (\ref{eq:E5}) and (\ref{eq:E4}) it follows
\begin{equation*}
\frac{c_{3}}{c_{5}}  =  \frac{((1-2t_{5}-2t_{1}+t_{5}^{2}+2t_{5}t_{1}+2t_{1}t_{3})S_{3}-t_{3}^{2})S_{1}t_{5}}{((1-2t_{3}-2t_{5}+t_{3}^{2}+2t_{3}t_{5}+2t_{5}t_{1})S_{1}-t_{1}^{2})s_{3}s_{5}t_{3}}.
\end{equation*}
 Using (\ref{eq:S3}) and (\ref{eq:S1}) we get
\begin{align}
&S_{3}+s_{3}s_{5}\frac{c_{3}}{c_{5}}=\nonumber\\
&-\frac{t_{3}S_{5}}
{\left(2(t_{5}-t_{3})(1-t_{1})S_{5}-t_{5}^{2}\right)\left((2t_{5}t_{1}-2t_{1}+2t_{1}t_{3}+t_{3}^{2}-2t_{3}+1)S_{5}-t_{5}^{2}\right)}\nonumber\\
&\times\left((t_{3}-2t_{3}^{2}-2t_{1}t_{3}+2t_{1}t_{3}^{2}+t_{3}^{3}+6t_{1}t_{3}t_{5}-2t_{1}t_{5}+t_{5}^{3}-2t_{3}t_{5}+t_{5})S_{5} \right.\nonumber\\
&\qquad \left.-t_{3}t_{5}^{2}-t_{5}^{3}\right).\label{eq:S5eq1left}
\end{align} 
Substituting (\ref{eq:S5eq1left}) into (\ref{eq:S5eq2}) we obtain
\begin{align*}
&\frac{6S_{5}(t_{3}+t_{5}-1)}{2t_{5}t_{1}-2t_{1}+2t_{1}t_{3}+t_{3}^{2}-2t_{3}+1)S_{5}-t_{5}^{2}}\\
&\times \left(t_{3}-2t_{3}^{2}-2t_{1}t_{3}+2t_{1}t_{3}^{2}+t_{3}^{3}+6t_{1}t_{3}t_{5}-2t_{1}t_{5}+t_{5}^{3}-2t_{3}t_{5}+t_{5})S_{5}\right.\\
&\qquad \left.-t_{3}t_{5}^{2}-t_{5}^{3}\right)\\
&=(1-t_{1})(t_{1}-5+6t_{3})S_{5}+3t_{5}^{2}.
\end{align*}
From this equality we can derive
\begin{align*}
&(5-72t_{1}t_{3}t_{5}-22t_{3}+35t_{3}^{2}-16t_{1}+52t_{1}t_{3}-54t_{1}t_{3}^{2}+22t_{1}t_{5}+24t_{3}t_{5}\\
&-12t_{3}t_{5}^{2}-6t_{5}-24t_{3}^{3}-6t_{5}^{3}    
+6t_{5}^{2}+13t_{1}^{2}t_{3}^{2}-24t_{3}^{2}t_{5}+6t_{5}t_{3}^{3}+18t_{1}t_{3}^{3}\\
&-12t_{1}^{2}t_{5}-26t_{1}^{2}t_{3}+2t_{1}^{3}t_{5}+2t_{1}^{3}t_{3}
+6t_{3}t_{5}^{3}+13t_{1}^{2}-12t_{1}t_{5}^{2}+6t_{5}^{4}-2t_{1}^{3}\\
&+6t_{3}^{4}+36t_{1}t_{3}t_{5}^{2}+48t_{1}t_{3}^{2}t_{5}+12t_{1}^{2}t_{3}t_{5})S_{5}^{2}\\
&+(-12t_{1}t_{3}t_{5}^{2}-6t_{5}^{3}t_{1}+12t_{1}t_{5}^{2}
-9t_{3}^{2}t_{5}^{2}-12t_{3}t_{5}^{3}+18t_{3}t_{5}^{2}
+6t_{5}^{3}-t_{1}^{2}t_{5}^{2}\\
&-8t_{5}^{2}-6t_{5}^{4})S_{5}+3t_{5}^{4}=0.
\end{align*}
Swapping the indices $1$ and $3$ we obtain
\begin{align*}
&(5-72t_{1}t_{3}t_{5}-22t_{1}+35t_{1}^{2}-16t_{3}+52t_{1}t_{3}-54t_{3}t_{1}^{2}+22t_{3}t_{5}+24t_{1}t_{5}\\
&-12t_{1}t_{5}^{2}-6t_{5}-24t_{1}^{3}-6t_{5}^{3}    
+6t_{5}^{2}+13t_{1}^{2}t_{3}^{2}-24t_{1}^{2}t_{5}+6t_{5}t_{1}^{3}+18t_{3}t_{1}^{3}\\
&-12t_{3}^{2}t_{5}-26t_{3}^{2}t_{1}+2t_{3}^{3}t_{5}+2t_{3}^{3}t_{1}
+6t_{1}t_{5}^{3}+13t_{3}^{2}-12t_{3}t_{5}^{2}+6t_{5}^{4}-2t_{3}^{3}\\
&+6t_{1}^{4}+36t_{1}t_{3}t_{5}^{2}+48t_{3}t_{1}^{2}t_{5}+12t_{3}^{2}t_{1}t_{5})S_{5}^{2}\\
&+(-12t_{1}t_{3}t_{5}^{2}-6t_{5}^{3}t_{3}+12t_{3}t_{5}^{2}
-9t_{1}^{2}t_{5}^{2}-12t_{1}t_{5}^{3}+18t_{1}t_{5}^{2}
+6t_{5}^{3}-t_{3}^{2}t_{5}^{2}\\
&-8t_{5}^{2}-6t_{5}^{4})S_{5}+3t_{5}^{4}=0.
\end{align*}
Since the right hand sides are equal in the above two equations and
$S_{5}\neq0$ we have
\begin{align}
&(t_{1}-t_{3})\left((3t_{1}^{3}+2t_{1}^{2}t_{5}+11t_{1}^{2}t_{3}-11t_{1}^{2}+11t_{1}+20t_{1}t_{3}t_{5}-6t_{1}t_{5}-25t_{1}t_{3}\right.\nonumber\\
&\qquad +11t_{1}t_{3}^{2}-3+2t_{3}^{2}t_{5}+3t_{5}^{3}+t_{5}-11t_{3}^{2}+3t_{3}^{3}+11t_{3}-6t_{3}t_{5})S_{5}\nonumber\\
&\qquad\left.-t_{5}^{2}(4t_{1}-3+3t_{5}+4t_{3})\right)=0.\label{eq:S5eq3}
\end{align}
Swapping the indices $3$ and $5$ we obtain
\begin{align}
&(t_{1}-t_{5})\left((3t_{1}^{3}+2t_{1}^{2}t_{3}+11t_{1}^{2}t_{5}-11t_{1}^{2}+11t_{1}+20t_{1}t_{5}t_{3}-6t_{1}t_{3}-25t_{1}t_{5}\right.\nonumber\\
&\qquad +11t_{1}t_{5}^{2}-3+2t_{5}^{2}t_{3}+3t_{3}^{3}+t_{3}-11t_{5}^{2}+3t_{5}^{3}+11t_{5}-6t_{5}t_{3})S_{3}\nonumber\\
&\qquad\left.-t_{3}^{2}(4t_{1}-3+3t_{3}+4t_{5})\right)=0. \label{eq:S3eq2}    
\end{align}
For later use (the case $t_{1}\neq t_{3}$) we note here that swapping
the indices $1$ and $5$ in (\ref{eq:S5eq3}) we obtain
\begin{align}
&(t_{5}-t_{3})\left((3t_{5}^{3}+2t_{5}^{2}t_{1}+11t_{5}^{2}t_{3}-11t_{5}^{2}+11t_{5}+20t_{5}t_{3}t_{1}-6t_{5}t_{1}-25t_{5}t_{3}\right.\nonumber\\
&\qquad +11t_{5}t_{3}^{2}-3+2t_{3}^{2}t_{1}+3t_{1}^{3}+t_{1}-11t_{3}^{2}+3t_{3}^{3}+11t_{3}-6t_{3}t_{1})S_{1}\nonumber\\
&\qquad\left. -t_{1}^{2}(4t_{5}-3+3t_{1}+4t_{3})\right)=0. \label{eq:S1eq2}
\end{align}

Assume that $t_{1}=t_{3}$. We investigate two cases.
Case A) $t_{1}\neq t_{5}$. \\
Then from (\ref{eq:S3eq2}) we obtain
\begin{align}
&\left(-17t_{1}^{2}+31t_{1}^{2}t_{5}+8t_{1}^{3}-31t_{1}t_{5}+13t_{1}t_{5}^{2}
+12t_{1}-11t_{5}^{2}+3t_{5}^{3}\right.\nonumber\\
&\left.-3+11t_{5}\right)S_{1} =t_{1}^{2}(7t_{1}+4t_{5}-3).\label{eq:S3eq3}
\end{align}
Substituting (\ref{eq:S1}) into (\ref{eq:S3eq3}) and using $t_{3}=t_{1}$,
$t_{1}\neq0$ we obtain
\begin{align}
&\left(-3-37t_{1}^{2}+25t_{1}^{2}t_{5}+22t_{1}^{3}-19t_{1}t_{5}+5t_{1}t_{5}^{2}+18t_{1}-3t_{5}^{2}
 \right. \nonumber\\
 &\left. +3t_{5}^{3}+5t_{5}\right)S_{5}=t_{5}^{2}(7t_{1}+4t_{5}-3).\label{eq:S5t1t3}
\end{align}
If the coefficient of $S_{5}$ is $0$ in (\ref{eq:S5t1t3}) then
we have
\begin{equation*}
-3-37t_{1}^{2}+25t_{1}^{2}t_{5}+22t_{1}^{3}-19t_{1}t_{5}+5t_{1}t_{5}^{2}+18t_{1}-3t_{5}^{2}+3t_{5}^{3}+5t_{5}=0,    
\end{equation*}
and
\begin{equation*}
7t_{1}+4t_{5}-3=0.   
\end{equation*}
From the last equality it follows
\begin{equation}
t_{5}=\frac{3-7t_{1}}{4}.\label{eq:t5}
\end{equation}
 Substituting it into the previous equation we obtain
\begin{equation*}
    (11t_{1}-3)(131t_{1}^{2}-42t_{1}+7)=0.
\end{equation*}
Since $131t_{1}^{2}-42t_{1}+7>0$, the only possibility is
\begin{equation*}
    t_{1}=\frac{3}{11}.
\end{equation*}
From (\ref{eq:t5}) we have
\begin{equation*}
    t_{5}=\frac{3}{11},
\end{equation*}
which contradicts the assumption $t_{1}\neq t_{5}$. Hence the coefficient
of $S_{5}$ cannot be $0$, and (\ref{eq:S5t1t3}) gives
\begin{align}
&S_{5}=\nonumber\\
&\frac{t_{5}^{2}(7t_{1}+4t_{5}-3)}{-3-37t_{1}^{2}+25t_{1}^{2}t_{5}+22t_{1}^{3}-19t_{1}t_{5}+5t_{1}t_{5}^{2}+18t_{1}-3t_{5}^{2}+3t_{5}^{3}+5t_{5}}.\label{eq:S5t1eqt3}
\end{align}
Using (\ref{eq:S5t1eqt3}) (and $t_{1}=t_{3}$, of course) we get
from (\ref{eq:S1}), (\ref{eq:C1}) and (\ref{eq:C5})
\begin{align}
&S_{1}=\nonumber\\
&\frac{t_{1}^{2}(7t_{1}+4t_{5}-3)}{-17t_{1}^{2}+31t_{1}^{2}t_{5}+8t_{1}^{3}-31t_{1}t_{5}+13t_{1}t_{5}^{2}+12t_{1}-11t_{5}^{2}+3t_{5}^{3}-3+11t_{5}},\label{eq:S1t1eqt3}
\end{align}
\begin{align}
&C_{1}=\nonumber\\
&\frac{(t_{1}-t_{5})(6t_{1}^{2}-3t_{1}-3t_{1}t_{5}+1-t_{5}^{2})}{2(-17t_{1}^{2}+31t_{1}^{2}t_{5}+8t_{1}^{3}-31t_{1}t_{5}+13t_{1}t_{5}^{2}+12t_{1}-11t_{5}^{2}+3t_{5}^{3}-3+11t_{5})},\label{eq:C1t1t3}
\end{align}
\begin{align}
&C_{5}=\nonumber\\
&\frac{(t_{1}^{2}-1-3t_{5}^{2}+3t_{5})}{2(6t_{1}^{2}-3t_{1}-3t_{1}t_{5}+1-t_{5}^{2})}\times\nonumber\\
&\frac{(t_{5}-t_{1}+3t_{1}t_{5}+t_{1}^{3}-t_{1}^{2}t_{5}-3t_{1}t_{5}^{2}-3t_{5}^{2}+3t_{5}^{3})}{(-3-37t_{1}^{2}+25t_{1}^{2}t_{5}+22t_{1}^{3}-19t_{1}t_{5}+5t_{1}t_{5}^{2}+18t_{1}-3t_{5}^{2}+3t_{5}^{3}+5t_{5})}.\label{eq:C5t1t3}
\end{align}
Since
\begin{eqnarray*}
S_{1}+C_{1}-1 & = & 0,\\
S_{5}+C_{5}-1 & = & 0,
\end{eqnarray*}
we obtain
\begin{align*}
q_{1}:&=4t_{1}^{3}-63t_{1}^{2}t_{5}+25t_{1}^{2}-23t_{1}-23t_{5}+65t_{1}t_{5}-24t_{1}t_{5}^{2}-5t_{5}^{3}+22t_{5}^{2}+6\\
&=0,
\end{align*}
\begin{align*}
q_{5}:&=6-11t_{5}-53t_{1}-484t_{1}^{3}+218t_{1}^{2}-3t_{5}^{3}-11t_{5}^{5}+18t_{5}^{4}-263t_{1}^{5}+576t_{1}^{4}\\    
&+212t_{5}^{2}t_{1}^{3}-114t_{1}^{2}t_{5}+44t_{1}t_{5}+85t_{1}t_{5}^{2}
-206t_{5}^{2}t_{1}^{2}+56t_{5}^{3}t_{1}^{2}-t_{5}^{4}t_{1}-62t_{5}^{3}t_{1}\\
&+162t_{1}^{3}t_{5}-169t_{1}^{4}t_{5}\\
&=0.
\end{align*}
To eliminate $t_{5}$ from the system we calculate the resultant
of $q_{1}$ and $q_{5}$ considering $t_{5}$ as the variable of these
two polynomials, and setting the resultant equal to zero.
\begin{eqnarray*}
R(q_{1},q_{5},t_{5}) & = & 48(38t_{1}-7)(61t_{1}^{2}-98t_{1}+49)(11t_{1}-3)^{2}\\
& &\times (131t_{1}^{2}-42t_{1}+7)^{2}(2t_{1}-1)^{6}\\
 & = & 0.
\end{eqnarray*}
 Similarly,
 \begin{eqnarray*}
R(q_{1},q_{5},t_{1}) & = & -48(38t_{5}-23)(61t_{5}^{2}+144t_{5}+11)(11t_{5}-3)^{2}\\
& &\times (131t_{5}^{2}-123t_{5}+40)^{2}(2t_{5}-1)^{6}\\
 & = & 0.
\end{eqnarray*}
Since $61t_{1}^{2}-98t_{1}+49,\,131t_{1}^{2}-42t_{1}+7,\,61t_{5}^{2}+144t_{5}+11,\,131t_{5}^{2}-123t_{5}+40>0$
we obtain the following system of equations 
\begin{equation*}
    (38t_{1}-7)(11t_{1}-3)(2t_{1}-1)=0,
\end{equation*} 
\begin{equation*}
    (38t_{5}-23)(11t_{5}-3)(2t_{5}-1)=0.
\end{equation*} 
By our assumptions $1>t_{2}+t_{3}=t_{1}+t_{3}=2t_{1}$ and $t_{1}\neq t_{5}$.
Thus we have to consider only the roots
\begin{equation*}
    (t_{1},t_{5})=\left(\frac{7}{38},\frac{23}{38}\right),\,\left(\frac{7}{38},\frac{3}{11}\right),\,\left(\frac{7}{38},\frac{1}{2}\right),\,\left(\frac{3}{11},\frac{23}{38}\right),\,\left(\frac{3}{11},\frac{1}{2}\right).
\end{equation*}
Now calculating the both sides in (\ref{eq:S1t1eqt3}) we see that
none of them is a solution of our original problem. ($\sin^{2}(7\pi/38)\neq147/211,\,5929/46828,\,539/4283$,
and $\sin^{2}(3\pi/11)\neq17328/7199,\,144/229$.)\\
Case B) $t_{1}=t_{5}$.\\
Substituting $t_{1}=t_{3}$ and $t_{1}=t_{5}$ into (\ref{eq:Eeq1})
and using $C_{1}=1-S_{1}$, it follows
\begin{equation}
S_{1}=\frac{3t_{1}^{2}}{7t_{1}^{2}-4t_{1}+1}.\label{eq:S1first}
\end{equation}
(Here $7t_{1}^{2}-4t_{1}+1>0$.) After that, substituting them into
(\ref{eq:Eeq7}) and using $C_{1}=1-S_{1}$, it follows
\begin{equation}
S_{1}=\frac{3t_{1}^{2}}{13t_{1}^{2}-9t_{1}+2}.\label{eq:S2second}
\end{equation}
 (Here $13t_{1}^{2}-9t_{1}+2>0$.) From (\ref{eq:S1first}) and (\ref{eq:S2second})
we obtain
\begin{equation*}
    t_{1}=\frac{1}{3}
\end{equation*}
or
\begin{equation*}
    t_{1}=\frac{1}{2}.
\end{equation*}
The last one is impossible, because $1>t_{2}+t_{3}=2t_{1}=1$. Substituting
the values $t_{i}=\frac{1}{3}$ ($i=1,\ldots,6$) into $A(\alpha,\beta)$
and introducing the new variables $\alpha=3x$ and $\beta=3y$ we
get 
\begin{align*}
&A(\alpha,\beta)=\\
&\sin^{2}(x)\sin^{2}(3x+3y)\cos^{2}\left(\frac{\pi}{6}+x\right)\cos^{2}\left(\frac{\pi}{6}+y\right)\\
&+\sin^{2}(3x)\cos^{2}\left(\frac{\pi}{6}+x+y\right)\sin^{2}(x+y)\cos^{2}\left(\frac{\pi}{6}+y\right)\\
&-2\sin(x)\sin(3x)\cos\left(\frac{\pi}{6}+x+y\right)\cos(y)\sin(x+y)\sin(3x+3y)\\
&\quad\times \cos\left(\frac{\pi}{6}+x\right)\cos^{2}\left(\frac{\pi}{6}+y\right)\\
&-\sin^{2}(3x)\sin^{2}(y)\sin^{2}(x+y)\cos^{2}\left(\frac{\pi}{6}+y\right)\\
&-\sin^{2}(3y)\sin^{2}(x)\sin^{2}(x+y)\cos^{2}\left(\frac{\pi}{6}+x\right)\\
&+2\sin(3x)\sin(y)\sin(3y)\sin(x)\sin\left(\frac{\pi}{6}+x+y\right)\sin^{2}(x+y)\\
&\quad\times\cos\left(\frac{\pi}{6}+x\right)\cos\left(\frac{\pi}{6}+y\right).
\end{align*}
Using the addition formulas of $\sin$ and $\cos$ functions 
 a very tedious calculation gives
\begin{equation*}
    A(\alpha,\beta)=0.
\end{equation*}
It means, if 
\begin{equation*}
    (t_{1},t_{2},t_{3},t_{4},t_{5},t_{6})=\left(\frac{1}{3},\frac{1}{3},\frac{1}{3},\frac{1}{3},\frac{1}{3},\frac{1}{3}\right),
\end{equation*}
then the triangle $GIJ_{\Delta}$ is equilateral. 

Now assume that $t_{1}\neq t_{3}$. 
In fact we may assume that $t_{1},t_{3},t_{5}$ are all different,
because the case $t_{1}=t_{3}$ can be considered as ''at least two
of the $t_{i}$'s are equal''. Then (\ref{eq:S5eq3}), (\ref{eq:S3eq2})
and (\ref{eq:S1eq2}) imply 
\begin{align}
&(3t_{1}^{3}+2t_{1}^{2}t_{5}+11t_{1}^{2}t_{3}-11t_{1}^{2}+11t_{1}+20t_{1}t_{3}t_{5}-6t_{1}t_{5}-25t_{1}t_{3}\nonumber\\
&+11t_{1}t_{3}^{2}-3+2t_{3}^{2}t_{5}+3t_{5}^{3}+t_{5}-11t_{3}^{2}+3t_{3}^{3}+11t_{3}-6t_{3}t_{5})S_{5}\nonumber\\
&=  t_{5}^{2}(4t_{1}-3+3t_{5}+4t_{3}),\label{eq:S5eqlast}
\end{align}
\begin{align}
&(3t_{1}^{3}+2t_{1}^{2}t_{3}+11t_{1}^{2}t_{5}-11t_{1}^{2}+11t_{1}+20t_{1}t_{5}t_{3}-6t_{1}t_{3}-25t_{1}t_{5}\nonumber\\
&+11t_{1}t_{5}^{2}-3+2t_{5}^{2}t_{3}+3t_{3}^{3}+t_{3}-11t_{5}^{2}+3t_{5}^{3}+11t_{5}-6t_{5}t_{3})S_{3}\nonumber\\
&= t_{3}^{2}(4t_{1}-3+3t_{3}+4t_{5}), \label{eq:S3eqlast}
\end{align}
\begin{align}
&(3t_{5}^{3}+2t_{5}^{2}t_{1}+11t_{5}^{2}t_{3}-11t_{5}^{2}+11t_{5}+20t_{5}t_{3}t_{1}-6t_{5}t_{1}-25t_{5}t_{3}\nonumber\\
&+11t_{5}t_{3}^{2}-3+2t_{3}^{2}t_{1}+3t_{1}^{3}+t_{1}-11t_{3}^{2}+3t_{3}^{3}+11t_{3}-6t_{3}t_{1})S_{1} \nonumber\\
&= t_{1}^{2}(4t_{5}-3+3t_{1}+4t_{3}).\label{eq:S1eqlast}
\end{align}
At this point we need some discussions. Assume that the coefficients
of $S_{i}$'s are equal to $0$. Then we obtain
\begin{eqnarray*}
4t_{1}-3+3t_{5}+4t_{3} & = & 0,\\
4t_{1}-3+3t_{3}+4t_{5} & = & 0,\\
4t_{5}-3+3t_{1}+4t_{3} & = & 0,
\end{eqnarray*}
whose solution is
\begin{equation*}
    t_{1},t_{3},t_{5}=\frac{3}{11},
\end{equation*}
which contradicts the assumption $t_{1}\neq t_{3}$. Thus one of the
coefficients, say that of $S_{5}$, cannot be $0$, and (\ref{eq:S5eqlast})
yields
\begin{align}
S_{5} &=  t_{5}^{2}(4t_{1}-3+3t_{5}+4t_{3})/(3t_{1}^{3}+2t_{1}^{2}t_{5}+11t_{1}^{2}t_{3}-11t_{1}^{2}+11t_{1}\nonumber\\
 &\qquad +20t_{1}t_{3}t_{5}-6t_{1}t_{5}-25t_{1}t_{3}+11t_{1}t_{3}^{2}-3+2t_{3}^{2}t_{5}+3t_{5}^{3}+t_{5}\nonumber\\
 &\qquad -11t_{3}^{2}+3t_{3}^{3}+11t_{3}-6t_{3}t_{5}).\label{eq:s5square}
\end{align}
If the coefficients of $S_{3}$ and $S_{1}$ are $0$, then
\begin{eqnarray*}
4t_{1}-3+3t_{3}+4t_{5} & = & 0,\\
4t_{5}-3+3t_{1}+4t_{3} & = & 0.
\end{eqnarray*}
 However this system gives $t_{1}=t_{3}$, a contradiction. Hence
without loss of generality we may assume that the coefficient of $S_{3}$
is not $0$, and (\ref{eq:S3eqlast}) gives 
\begin{align}
S_{3} & = t_{3}^{2}(4t_{1}-3+3t_{3}+4t_{5})/(3t_{1}^{3}+2t_{1}^{2}t_{3}+11t_{1}^{2}t_{5}-11t_{1}^{2}\nonumber\\
&\qquad +11t_{1}+20t_{1}t_{5}t_{3}-6t_{1}t_{3} -25t_{1}t_{5}+11t_{1}t_{5}^{2}-3\nonumber\\
&\qquad +2t_{5}^{2}t_{3}+3t_{3}^{3}+t_{3}-11t_{5}^{2}
+3t_{5}^{3}+11t_{5}-6t_{5}t_{3}).\label{eq:s3squareV1} 
\end{align}
Now we distinguish two cases.

Case A) The coefficient of $S_{1}$ is $0$.

Then we obtain 
\begin{eqnarray}
3t_{5}^{3}+2t_{5}^{2}t_{1}+11t_{5}^{2}t_{3}-11t_{5}^{2}+11t_{5}+20t_{5}t_{3}t_{1}-6t_{5}t_{1}-25t_{5}t_{3}\nonumber \\
+11t_{5}t_{3}^{2}-3+2t_{3}^{2}t_{1}+3t_{1}^{3}+t_{1}-11t_{3}^{2}+3t_{3}^{3}+11t_{3}-6t_{3}t_{1} & = & 0,\label{eq:CaseAt1t3t5}
\end{eqnarray}
and
\begin{equation*}
4t_{5}-3+3t_{1}+4t_{3}=0.
\end{equation*}
From the last one we have
\begin{equation}
t_{5}=\frac{3-3t_{1}-4t_{3}}{4}.\label{eq:t5eq}
\end{equation}
 Substituting (\ref{eq:t5eq}) into (\ref{eq:CaseAt1t3t5}) we get
\begin{equation}
-\frac{131}{64}t_{1}-\frac{9}{8}t_{3}-\frac{9}{64}t_{1}^{2}+\frac{183}{64}t_{1}^{3}+\frac{3}{2}t_{3}^{2}-\frac{87}{8}t_{1}^{2}t_{3}+12t_{1}t_{3}-\frac{29}{2}t_{1}t_{3}^{2}+\frac{21}{64}=0.\label{eq:CaseAt1t3first}
\end{equation}
Furthermore, substituting (\ref{eq:t5eq}) into (\ref{eq:s5square})
and (\ref{eq:s3squareV1}) we obtain
\begin{align}
S_{5} &=(-3+3t_{1}+4t_{3})^{2}(7t_{1}-3+4t_{3})/
(125t_{1}+28t_{3}-77t_{1}^{2}+15t_{1}^{3}\nonumber\\
&\qquad +208t_{3}^{2}-128t_{3}^{3}-708t_{1}^{2}t_{3}+680t_{1}t_{3}
 -1104t_{1}t_{3}^{2}-63),\label{eq:CaseAS5}
\end{align}
\begin{align}
S_{3} &=-64t_{3}^{2}(t_{1}-t_{3})/
 (79t_{1}+124t_{3}-79t_{1}^{2}+21t_{1}^{3}+80t_{3}^{2}\nonumber\\
&\qquad -128t_{3}^{3}+732t_{1}^{2}t_{3}-856t_{1}t_{3}+816t_{1}t_{3}^{2}-21).\label{eq:CaseAS3}
\end{align}
Substituting (\ref{eq:t5eq}) and (\ref{eq:CaseAS5}) into (\ref{eq:S3})
we get
\begin{align}
S_{3} &=-16t_{3}^{2}(7t_{1}-3+4t_{3})/
(-331t_{1}^{2}+187t_{1}-1512t_{1}t_{3}-9+260t_{3}\nonumber\\
&\qquad +153t_{1}^{3}+1252t_{1}^{2}t_{3}+1360t_{1}t_{3}^{2}-464t_{3}^{2}+128t_{3}^{3}).\label{eq:AS3v2}
\end{align}
From (\ref{eq:CaseAS3}) and (\ref{eq:AS3v2}) it follows
\begin{align}
&-16t_{3}^{2}(-3+3t_{1}+8t_{3})
(155t_{1}^{3}-684t_{1}^{2}t_{3}-81t_{1}^{2}-95t_{1}\nonumber\\
&\qquad +792t_{1}t_{3}-912t_{1}t_{3}^{2}-108t_{3}+21+144t_{3}^{2})  =  0.\label{eq:CaseAt1t3second}
\end{align}
Here there are two possibilities.

(i) $-3+3t_{1}+8t_{3}=0$.\\
Then
\begin{equation}
t_{3}=\frac{3-3t_{1}}{8}.\label{eq:CaseAit3}
\end{equation}
Substituting it into (\ref{eq:CaseAt1t3first}) we obtain
\begin{equation}
(11t_{1}-3)(57t_{1}^{2}-36t_{1}-5)=0.
\end{equation}
 Its roots are
\[
t_{1}=\frac{3}{11},\frac{6}{19}+\frac{\sqrt{609}}{57},\frac{6}{19}-\frac{\sqrt{609}}{57}.
\]
 Here the third root is negative. Using (\ref{eq:CaseAit3}) we obtain
\begin{equation*}
t_{3}=\frac{3}{11},\frac{39}{152}-\frac{\sqrt{609}}{152}.
\end{equation*}
The first is impossible because $t_{1}\neq t_{3}$, so the only one
possibility is
\begin{equation*}
(t_{1},t_{3})=\left(\frac{6}{19}+\frac{\sqrt{609}}{57},\frac{39}{152}-\frac{\sqrt{609}}{152}\right).
\end{equation*}
Substituting these values into (\ref{eq:t5eq}) we get
\begin{equation*}
t_{5}=\frac{39}{152}-\frac{\sqrt{609}}{152},
\end{equation*}
which is impossible because $t_{3}\neq t_{5}.$

(ii) $-3+3t_{1}+8t_{3}\neq0$.\\
Then we have
\begin{align*}
p_{1}:=&-\frac{131}{64}t_{1}-\frac{9}{8}t_{3}-\frac{9}{64}t_{1}^{2}+\frac{183}{64}t_{1}^{3}+\frac{3}{2}t_{3}^{2}-\frac{87}{8}t_{1}^{2}t_{3}\\
&+12t_{1}t_{3}-\frac{29}{2}t_{1}t_{3}^{2}+\frac{21}{64}=0,
\end{align*}
\begin{align*}
p_{3}:=&155t_{1}^{3}-684t_{1}^{2}t_{3}-81t_{1}^{2}-95t_{1}+792t_{1}t_{3}\\
&-912t_{1}t_{3}^{2}-108t_{3}+21+144t_{3}^{2}=0.
\end{align*}
Calculating the resultants we obtain
\begin{equation*}
R(p_{1},p_{3},t_{3})=\frac{1}{16}(11t_{1}-3)^{2}(t_{1}+3)^{2}(131t_{1}^{2}-42t_{1}+7)^{2}=0,
\end{equation*}
\begin{align*}
R(p_{1},p_{3},t_{1})=&-9(5t_{3}^{2}-15t_{3}-8)(131t_{3}^{2}-42t_{3}+7)\\
& \cdot (131t_{3}^{2}-123t_{3}+40)(11t_{3}-3)=0.
\end{align*}
Since $t_{1}+3,131t_{1}^{2}-42t_{1}+7>0$ the only one possibility
is $t_{1}=3/11$. Similarly, since $131t_{3}^{2}-42t_{3}+7,131t_{3}^{2}-123t_{3}+40>0$,
and $5t_{3}^{2}-15t_{3}-8<0$ the only one possibility is $t_{3}=3/11$,
but it is impossible because $t_{1}\neq t_{3}$.

Case B) The coefficient of $S_{1}$ is not $0$.\\
Then 
\begin{align}
S_{1} = & t_{1}^{2}(4t_{5}-3+3t_{1}+4t_{3})/(3t_{5}^{3}+2t_{5}^{2}t_{1}+11t_{5}^{2}t_{3}-11t_{5}^{2}+11t_{5}\nonumber\\
&\quad +20t_{5}t_{3}t_{1}-6t_{5}t_{1}-25t_{5}t_{3}
 +11t_{5}t_{3}^{2}-3+2t_{3}^{2}t_{1}+3t_{1}^{3}\nonumber\\
 &\quad +t_{1}-11t_{3}^{2}+3t_{3}^{3}+11t_{3}-6t_{3}t_{1}).\label{eq:s1squareV2}
\end{align}
Using (\ref{eq:s5square}) we get from (\ref{eq:S1}) and (\ref{eq:S3})
\begin{align}
S_{1} = & t_{1}^{2}(4t_{1}-3+3t_{5}+4t_{3})/(-3+22t_{1}t_{3}t_{5}+11t_{3}-11t_{3}^{2}+5t_{1}\nonumber\\
&\quad -11t_{1}t_{3}+3t_{1}t_{3}^{2}-8t_{1}t_{5}
-20t_{3}t_{5}+6t_{3}t_{5}^{2}+7t_{5}+3t_{3}^{3}\nonumber\\
&\quad +3t_{5}^{3}-6t_{5}^{2}+10t_{3}^{2}t_{5}+2t_{1}^{2}t_{5}+3t_{1}^{2}t_{3}-3t_{1}^{2}+3t_{1}^{3}),\label{eq:s1square}
\end{align}
\begin{align}
S_{3} = & t_{3}^{2}(4t_{1}-3+3t_{5}+4t_{3})/(-3+22t_{1}t_{3}t_{5}+11t_{1}-11t_{1}^{2}+5t_{3}\nonumber\\
&\quad -11t_{1}t_{3}+3t_{3}t_{1}^{2}-8t_{3}t_{5}
-20t_{1}t_{5}+6t_{1}t_{5}^{2}+7t_{5}+3t_{1}^{3}\nonumber\\
&\quad +3t_{5}^{3}-6t_{5}^{2}+10t_{1}^{2}t_{5}+2t_{3}^{2}t_{5}+3t_{3}^{2}t_{1}-3t_{3}^{2}+3t_{3}^{3}).\label{eq:s3square}
\end{align}
Since the left hand sides are the same, we obtain from (\ref{eq:s1squareV2})
and (\ref{eq:s1square})
\begin{align*}
t_{1}^{2}(t_{1}-t_{5})(& 3t_{1}^{3}+9t_{1}^{2}-6t_{1}^{2}t_{5}-9t_{1}^{2}t_{3}+27t_{5}t_{1}+30t_{1}t_{3}-20t_{1}-13t_{1}t_{3}^{2}\\
& -15t_{1}t_{3}t_{5}
-6t_{5}^{2}t_{1}+9-t_{3}^{3}-20t_{3}-9t_{3}t_{5}^{2}\\
&+9t_{5}^{2}+3t_{5}^{3}+12t_{3}^{2}+30t_{3}t_{5}-20t_{5}-13t_{3}^{2}t_{5})\\
& = 0.
\end{align*}
Similarly, we obtain from (\ref{eq:s3squareV1}) and (\ref{eq:s3square})
\begin{align*}
t_{3}^{2}(t_{3}-t_{5})(&-3t_{3}^{3}+9t_{1}t_{3}^{2}-9t_{3}^{2}+6t_{3}^{2}t_{5}-30t_{1}t_{3}-27t_{3}t_{5}+13t_{1}^{2}t_{3}\\
& +15t_{1}t_{3}t_{5}+20t_{3}
+6t_{3}t_{5}^{2}-9+t_{1}^{3}+20t_{1}-9t_{5}^{2}-3t_{5}^{3}\\
& -12t_{1}^{2}+9t_{5}^{2}t_{1}-30t_{5}t_{1}+13t_{1}^{2}t_{5}+20t_{5})\\
& = 0.
\end{align*}
Since $t_{1},t_{3},t_{1}-t_{5},t_{3}-t_{5}\neq0$ we obtain
\begin{align*}
& 3t_{1}^{3}+9t_{1}^{2}-6t_{1}^{2}t_{5}-9t_{1}^{2}t_{3}+27t_{5}t_{1}+30t_{1}t_{3}-20t_{1}-13t_{1}t_{3}^{2}-15t_{1}t_{3}t_{5}-6t_{5}^{2}t_{1}\\
& +9-t_{3}^{3}-20t_{3}-9t_{3}t_{5}^{2}
+9t_{5}^{2}+3t_{5}^{3}+12t_{3}^{2}+30t_{3}t_{5}-20t_{5}-13t_{3}^{2}t_{5}\\
& = 0,
\end{align*}
and
\begin{align*}
&-3t_{3}^{3}+9t_{1}t_{3}^{2}-9t_{3}^{2}+6t_{3}^{2}t_{5}-30t_{1}t_{3}-27t_{3}t_{5}+13t_{1}^{2}t_{3} +15t_{1}t_{3}t_{5}+20t_{3}\\
& +6t_{3}t_{5}^{2}-9+t_{1}^{3}+20t_{1}-9t_{5}^{2}-3t_{5}^{3}
 -12t_{1}^{2}+9t_{5}^{2}t_{1}-30t_{5}t_{1}+13t_{1}^{2}t_{5}+20t_{5}\\
& = 0.
\end{align*}
 Summing the left and right sides we have
\begin{equation*}
(t_{1}-t_{3})(t_{1}+t_{3}+t_{5})(4t_{1}-3+3t_{5}+4t_{3})=0.
\end{equation*}
Here $t_{1}-t_{3},t_{1}+t_{3}+t_{5}\neq0$, consequently
\begin{equation*}
4t_{1}-3+3t_{5}+4t_{3}=0.
\end{equation*}
Substituting it into (\ref{eq:S5eqlast}) we obtain that the coefficient
of $S_{5}$ must be $0$ which is a contradiction. Hence there is
no solution in Case B.
\end{proof}

\end{document}